\documentclass{notices}
\usepackage{amsfonts,amssymb,amsmath,amscd}

\title{The Mathematician Leads\\
  \large Building a Research Profile in the Age of Frontier AI}

\author{Helmut H. W. Hofer
  \affil{Helmut H. W. Hofer is Professor Emeritus in the School of
  Mathematics at the Institute for Advanced Study. His email address is
  \texttt{hofer@ias.edu}.}}

\date{}

\begin{document}

\maketitle

\section{The proposal}

In 2025 I tried a frontier AI model in my mathematical work. It helped with
routine tasks but not as a serious research collaborator. A newer model changed
my assessment. For several months I have used it in ongoing research: for
calculations, proof attempts, manuscript construction, literature questions,
changes of strategy, and the diagnosis of failed approaches.

I learned quickly that access to a strong model is not yet a research
collaboration. The model often knew the relevant definitions, theorems, and
standard constructions, but it could still lose the direction of an emerging
argument. Faced with uncertainty, it sometimes moved toward a familiar
neighboring problem, polished an incomplete idea, or invoked a standard
framework before checking whether the new mechanism fit it.

We began to record my corrections and to ask whether they revealed recurring
patterns. Those records developed into what I call a \emph{research profile}:
a researcher-controlled, versioned account of how a particular mathematician
and an AI system work together. It contains hypotheses about the collaboration,
the instructions currently being tested, the episodes that motivated them,
and the evidence that later supports or contradicts them. Nothing in the
profile is permanent. Every instruction has a kill switch: if it makes the
collaboration worse, it should be removed.

A profile is a falsifiable account of how a collaboration currently works
best, not a certificate and not a permanent portrait of the mathematician.
Its value must be tested on later work.

One recurring correction concerned large-step connections between distant
areas. The model's breadth helped it find relevant material, but also tempted
it toward nearby problems that were easier and more familiar. Several such
episodes produced this instruction:

\begin{quote}
\itshape Preserve the intended direction of a large-step idea before
decomposing it. Do not replace a difficult proposed relation with a familiar
neighboring problem merely because the latter is easier.
\end{quote}

For each instruction we record where it came from, what kinds of sessions it
governs, and whether it prevents the same drift in later work. The evidence and
revision history are as important as the instruction itself.

I want to know whether the corrections and judgments accumulated in a sustained
mathematical collaboration can become a persistent, testable profile, and
whether generic rules can be separated from practices reflecting one
mathematician's judgment. I describe one developing experiment and propose ways
to test it; I do not claim that the case has been established.

\section{How the experiment began}

Frontier AI models search the literature, calculate, produce code, and rewrite
arguments with remarkable speed. Speed is not research judgment. Public
discussion nevertheless treats success on one problem as evidence for
automation, or a conspicuous failure as evidence against it. Present
capabilities are uneven: a model may excel when the framework is visible and
lose direction while the formulation, relevant structures, and standards of
evidence are changing. The practical question is how mathematicians can use AI
while remaining the intellectual leaders of their work.

Before this collaboration began, I had spent about two years developing, in
the traditional ``brain-assisted'' way, the mathematical landscape for a
long-term project in symplectic and contact topology. The AI neither supplied
that landscape nor chose its direction. I first gave the model a typeset
manuscript proving one direction of a local version of the project. I asked for
an adversarial audit: search for
missing hypotheses, sign errors, hidden circularity, and analogies standing in
for arguments. The proof survived, and the model shortened several passages
without weakening them.

We then turned to the difficult converse direction: proposals, objections,
discarded approaches, reformulations, and useful partial structures. At this
stage our discussion of the local case seems to be converging toward a proof
whose formulation may help reach the global version, the project's goal. The
incomplete argument must still be judged as mathematics.

The model's breadth was useful because the project crossed established areas,
but the same breadth caused drift. The large-step instruction quoted above came
from this experience. A second recurring issue concerned agreement. When I
asked, ``Do you agree?'' I was requesting an independent check of signs,
dimensions, symmetries, domains, and logical implications. That correction
became another compact rule:

\begin{quote}
\itshape An invitation to agree is an invitation to audit.
\end{quote}

When the model proposed a relation between two areas, I also required it to
identify the objects being compared, the dictionary between them, the
structures preserved, the necessary hypotheses, and the likely point of
failure. A verbal resemblance was not yet a mathematical bridge.

Other instructions ask the model to build connected mathematical landscapes
rather than inventories of facts, distinguish analogies from established
results, record why attractive approaches failed, and treat contrarian
collaboration as a positive duty.

These rules were not written in advance. They grew out of mistakes and
corrections in ordinary work. This article arose in the same way. Codex
produced a first draft; I criticized its emphasis, scope, and presentation; and
we rewrote it repeatedly. That history illustrates the method, but it does not
prove that the method works.

\section{What a research profile contains}

A research profile is neither a curriculum vitae nor a personality test, and
it does not attempt to reproduce a mathematician inside a machine. In practice
we keep four records separate.

\begin{itemize}
\item The \emph{profile} contains current hypotheses about the mathematician,
  the AI system, and their interaction.
\item The \emph{manual} contains the instructions that are presently active.
\item The \emph{casebook} records the episodes, corrections, and later evidence
  supporting or opposing those instructions.
\item The \emph{project memory} records the mathematical state of a particular
  investigation: definitions, open implications, abandoned routes, notation,
  and manuscript versions.
\end{itemize}

One casebook entry records an incorrect symmetry reduction. The model
quotiented one visible factor by a circle action and obtained an attractive
familiar space, although the circle acted diagonally on two parts of the
configuration. My correction became a diagnostic question: before forming a
quotient, write the action on every variable and identify the full quotient
map. The manual keeps the instruction, the casebook the episode, and the
project memory the calculation. Later uses show whether the instruction helps
or merely repairs one mistake.

A substantial profile change should record:

\begin{enumerate}
\item the episode that motivated the change;
\item the diagnosed weakness or opportunity;
\item the proposed instruction;
\item the expected effect and the situations in which it applies;
\item evidence from later use; and
\item the condition under which the instruction should be reversed.
\end{enumerate}

A researcher or graduate student develops judgment in much the same way: not
only by collecting correct statements, but also by remembering which
attractive routes failed, which questions exposed the difficulty, and which
habits improved later work. I do not mean that the model thereby acquires
human understanding. The profile is an external record, and its value lies in
remaining visible, contestable, and under the mathematician's control.

Instructions come from different layers. Generic rules concern verification,
provenance, and uncertainty; individual rules concern standards of evidence or
ways of moving between local and global viewpoints; other rules repair a
current model or serve one project. The layers must remain distinct, and the
same parts need not govern discovery, proof verification, and exposition.

A mathematician need not have exactly one profile. Discovery, proof
verification, and exposition place different demands on the collaborator and
may activate different parts of the same profile.

\section{The profile changes the mathematician as well}

After we had worked together for a month or two, I asked the model what it had
learned about me as a mathematician. Its answer was offered for discussion,
not as a diagnosis. It observed that I use short prompts as probes in one
evolving argument, move repeatedly between local models and global structures,
and search for a decisive mechanism before accepting a large accumulation of
computations.

I printed the response and showed it to my wife, who has a background in
psychology and to whom I have been married for forty-three years. She laughed
and said that she could have told me the same things. The model had not
uncovered a hidden personality. It had reconstructed recognizable features of
my way of thinking from the structure of our conversations.

The episode made my own working habits explicit enough to examine. Do I choose
a local question because its answer would change the global picture, or merely
because it is accessible? Am I attracted to a construction because it is
structurally appropriate, or because I know it well? Which corrections concern
the model, and which reveal something about my own practice? A profile does not
answer these questions; it makes them harder to avoid.

Most changes grew out of ordinary research sessions and short discussions
afterward. Codex organized the record and drafted formulations; I accepted,
changed, or rejected them. The profile accumulated as a byproduct and then
began to influence the research in return.

I began to think that the effective unit of mathematical research may be
changing. As strong models become widely available, access alone will cease to
distinguish researchers. Two mathematicians using the same model may obtain
very different value from it. The relevant unit may increasingly be

\[
  \text{mathematician}+\text{research profile}+\text{frontier model}.
\]

The model can be replaced, but the profile cannot be assumed to transfer
unchanged. Generic and mathematician-specific layers may survive; repairs
aimed at one model may not. A profile should therefore be retested when the
model changes.

\section{Neighboring ideas and the missing layer}

Tao's survey of machine-assisted proof places current AI systems in a longer
history of mathematical computation \cite{Tao2025}. It treats machine learning,
proof assistants, and language models as tools with complementary strengths,
and envisages assistants that suggest ideas, filter errors, and perform routine
work while mathematicians concentrate on high-level concepts. Henkel similarly
argues for ``copilot, not pilot,'' combining strategic prompting with critical
verification and an experimental attitude toward models \cite{Henkel2025}.

The problem discussed here arises after that division of labor has been
accepted. Once a model enters a sustained and unresolved research program, can
corrections and judgments from the work improve later collaboration in a
controlled and testable way?

Two recent mathematical systems are relevant. The AI co-mathematician provides
a stateful workspace that manages uncertainty,
refines user intent, tracks failed hypotheses, and produces evolving
mathematical artifacts \cite{Zheng2026}. MathCoPilot keeps the mathematician in
control of an editable proof blueprint, combines proving strategies, and builds
a personal Lean-verified knowledge base \cite{Zhang2026}.

Neither system quite addresses the layer I have in mind. A profile can tell the
AI that a request for agreement requires adversarial checking; that exploratory
work should preserve a proposed relation instead of normalizing it into a
familiar framework; that a literature search should report conflicting
conventions rather than silently choose one; or that an exposition session may
compress a settled argument but must not settle an unresolved one. The issue
is when and how tools are used, not only which tools are available.

Personalization is also an active subject in AI research. Personalized Deep
Research incorporates user context into query development, retrieval, and
synthesis \cite{Li2026}. More directly, Huang, Du, and Lan distilled
natural-language skills from 206 sessions involving 13 developers and
evaluated them by simulated replay on held-out coding tasks \cite{Huang2026}.
The pooled generic skill performed best, but its advantage over no skill
narrowly missed the conventional significance threshold; developer-specific
gains were smaller and likewise inconclusive. This evidence is too thin for a
general conclusion. It does make one question testable: with more interaction
history, when do individual profiles outperform generic ones?

Open-ended mathematical research may behave differently. Coding tasks usually
have fixed objectives and tests; in research, the problem, formulation, and
standards of progress may change. One should therefore compare four
conditions: no persistent profile, a generic manual, the mathematician's own
profile, and a profile derived from another researcher. Testing must occur on
later work, not on the episodes from which a profile was constructed.

Rules for verification and provenance may transfer broadly. Whether the same
is true of question selection, mathematical taste, or movement between local
and global viewpoints is less clear. The issue becomes sharper with imagined
profiles of Jean Bourgain, Andrew Wiles, and Grigori Perelman. If the same
generic manual could replace all three without changing the questions pursued,
the structures emphasized, or the development of an argument, then the
individual layer would have little value. This should be tested rather than
assumed.

The hypothesis is modest but testable: some of the working judgment developed
in a mathematician--model collaboration can be recorded, revised, and tested
on later work. Generic and individual components should be tested separately.

\section{What profiles might be good for}

The first use of a profile is in a mathematician's own work with AI. Improvement
need not mean producing more pages or proving results faster. A useful rule may
prevent a repeated conceptual error, preserve the
direction of an argument, expose a hidden assumption earlier, improve the
choice of the next local question, or help recover from a failed path.

Career stage matters. An established researcher may use a profile to preserve
a distinctive direction while exploiting the model's breadth. A younger
mathematician may benefit from protocols for
checking claims and separating analogy from theorem, but may also be more
vulnerable to confident redirection. A borrowed profile must not substitute
another person's taste for the development of one's own judgment.

A complete profile may remain private while its owner shares an instruction for
adversarial checking, a sanitized failed analogy, or an evaluation problem.
Such components could help others without exposing unpublished mathematics or
pretending that one person's profile is universal.

For a group, a collaborative profile can supplement individual profiles by
specifying common language, division of work, standards of evidence,
permissions, and procedures for disagreement. A geometer, analyst, algebraist,
and formalizer need not be homogenized. The AI could use one perspective,
compare several, or return a disagreement to the humans. Mathematics often
depends on this productive tension, not on one averaged method.

Success must be judged differently in each setting: does an individual profile
reduce repeated corrections or expose a decisive question earlier; does a
younger researcher gain independent judgment rather than compliance; does a
group clarify disagreement without erasing distinct strengths? These measures
do not establish mathematical success, but they give claims that can be
examined rather than advertised.

\section{Profiles must not become cages}

A profile can fossilize present habits, teach the model to imitate style rather
than improve judgment, reinforce blind spots, or give unwarranted authority to
the template of a famous mathematician. Part of a good profile should instruct
the model to challenge the profile itself.

Using a corpus of 41.3 million English-language papers from six natural-science
disciplines published between 1980 and 2025, a language-model classifier
identified roughly 311,000 as AI-augmented. The observational study associated
AI use with greater individual output and citation impact but a narrower
collective range of topics \cite{Hao2026}. It establishes neither causation nor
anything directly about research profiles, but it raises a relevant question:
will personalization preserve distinctive structures of attention or narrow
them further?

Profiles also raise questions of autonomy and privacy. In a study of 668
ChatGPT users, fine-tuned classifiers inferred some personality traits from
conversation histories at rates above the reported baselines, with accuracy
depending strongly on the conversational setting \cite{Cogendez2026}. More
generally, personalized alignment can create risks of profiling, privacy
infringement, bias reinforcement, and manipulation \cite{Kirk2024}.
Unpublished mathematics may identify a researcher even after names are
removed. No serious program should promise confidentiality while silently
using every interaction for shared training.

Ongoing research should remain private by default. A researcher may separately
authorize an abstracted lesson, evaluation problem, correction pattern, or
published episode, and should see the precise object being shared. Academic
institutions should not demand private collaboration records or turn profiles
into employment metrics.

\section{A grassroots experiment}

I would begin modestly: mathematicians could construct small profiles and test
them in their ordinary work. No technical setup is required to begin;
provisional instructions can emerge from an ordinary conversation with the
model.

\begin{enumerate}
\item Use a frontier model in ordinary mathematical work and keep a short
  record of substantive corrections.
\item After several sessions, identify two or three recurring behaviors that
  help or hinder the collaboration.
\item Write a compact initial profile containing only instructions motivated
  by those observations.
\item Test the profile on new work rather than on the episodes from which it
  was derived.
\item Record where the intervention helped, had no effect, or produced a new
  failure.
\item Revise or remove the instruction accordingly.
\end{enumerate}

One experiment might ask whether ``an invitation to agree is an invitation to
audit'' produces substantive objections without empty contrarianism. Another
might test whether the large-step instruction prevents drift on new cross-field
problems without impeding ordinary calculations. A third could compare a
generic proof-checking manual with a mathematician-specific profile during
discovery and exposition. The useful questions are conditional: what helps,
for whom, at which stage, and at what cost?

The experiment can remain private. A broader community could share selected
profile components, sanitized cases, evaluation problems, and negative
findings. A profile that changes only tone, consumes too much time, reinforces
a blind spot, or fails after a model update is informative.

A small group of established researchers could also ask which components
survive model changes, which lessons transfer, and whether collaborative
profiles coordinate different perspectives without exposing unpublished
mathematics.

\section{The mathematician must lead}

AI is already entering mathematics through the daily work of individual
researchers. The relevant choice is not between mathematics with AI and
mathematics without it, but whether mathematicians develop the intellectual
means to govern the collaboration.

Our experiment remains unfinished, and I expect some current rules to
disappear. A living profile should develop through research. For now the
proposal is modest. Mathematicians who are interested should
construct small profiles, retain control of their research, compare what works,
and report failures as carefully as successes.

Models can of course generate ideas. The mathematician leads by remaining
responsible for the agenda, the standards of evidence, the decision to pursue
or abandon a direction, and the final claims. A profile may make that
collaboration more critical, personal, and mathematically effective. A broader
practice should grow from evidence rather than prescription.

\section*{Disclosure}

The collaboration described here was conducted with OpenAI Codex using
frontier models during 2026. Codex assisted with mathematical discussion,
organization, literature searches, reference auditing, and drafting. I
selected the direction, evaluated and revised the mathematical and
methodological claims, and take responsibility for the article.

\bibliography{Research_PROFILE_arxiv}

\end{document}